\documentclass[12pt]{amsart}
\usepackage{amsmath, amsthm, amscd, amsfonts}

\makeatletter \oddsidemargin.9375in \evensidemargin \oddsidemargin
\newtheorem{theorem}{Theorem}[section]
\newtheorem{lemma}[theorem]{Lemma}

\theoremstyle{definition}
\newtheorem{definition}[theorem]{Definition}

\theoremstyle{remark}

\numberwithin{equation}{section}

\begin{document}
\title[ probabilistic N-Banach spaces ] {A note on the probabilistic N-Banach spaces}
\author[POURMOSLEMI and SALIMI]{\textbf{ALIREZA POURMOSLEMI}\\Department of Mathematics, Payame Noor University, Iran.\\
pourmoslemy@yahoo.com
\\and\\
\textbf{MEHDI SALIMI}\\Department of Mathematics, Islamic Azad University, Toyserkan Branch, Iran.\\
msalimi1@yahoo.com\\}


\begin{abstract}
In this paper, we define probabilistic n-Banach spaces along with
some concepts in this field and study convergence in these spaces
by some lemmas and theorem.
\end{abstract} \maketitle

\section{Introduction and Preliminaries}
\par{In \cite{me} K. Menger introduced the notion of probabilistic metric
spaces. The idea of Menger was to use distribution function
instead of nonnegative real numbers as values of the metric. The
concept of probabilistic normed spaces were introduced by
\v{S}erstnev in \cite{ser}. New definition of probabilistic normed
spaces and some concepts in this field were studied by Alsina,
Schweizer, Sklar, Pourmoslemi and Salimi
\cite{al1,al,al2,pour1,pour2}. It corresponds to the situations
when we do not know exactly the distance between two points, we
know only probabilities of possible values of this distance.
 The probabilistic generalization of metric
spaces appears to be well adapted for the investigation of quantum
particle physics particularly in connections with both string and
$\varepsilon^\infty$ theory which were given and studied by El
Naschie \cite{{na1},{na2}}. A {\it distribution function (briefly a
~d.f.)} is a function $F$ from the extended real line
$\overline{\mathbb{R}}=[-\infty, +\infty]$ into the
 unit interval $I=[0, 1]$ that is nondecreasing and satisfies
$F(-\infty)=0,~ F(+\infty)=1.$ The set of all d.f.'s will be denoted
by $\Delta$ and the subset of those d.f.'s such that $F(0)=0$, will
be  denoted by $\Delta^+$ and $D^+\subseteq\Delta^+$ is defined as
follows:
\begin{equation*}
D^+=\{F\in\Delta^+:l^-F(+\infty)=1\},
\end{equation*}
where $l^-f(x)$ denotes the left limit of the function $f$ at the
point $x$. By setting $F\leq G$ whenever $F(x)\leq G(x)$ for all $x$
in $\mathbb{R}$, the maximal element for $\Delta^+$ in this order is
the d.f. given by
\[ \mathcal H_0(x) = \left\lbrace
  \begin{array}{c l}
    0, & \text{if $x\leq 0$},\\
    1, & \text{if $x>0$}.
  \end{array}
\right. \]

The space $\Delta$ can be metrized in several equivalent ways
\cite{{sc},{se},{si},{ty}} in such a manner that the metric topology
coincides with the topology of weak convergence for distribution
functions. Here, we assume that $\Delta$ is metrized by the Sibley
metric $d_s$. If $F$ and $G$ are d.f.'s and $h$ is in $]0, 1]$, let
$(F,G;h)$ denote the condition
\begin{eqnarray*}
F(x-h)-h\leq G(x)\leq F(x+h)+h,~~~~~~x\in
]-\frac{1}{h},\frac{1}{h}].
\end{eqnarray*}
Then the Sibley metric is defined by
\begin{eqnarray*}
d_s(F,G):= inf\{h\in ]0,1]:~ both~~~ (F,G;h)~~~and ~~~(G,F;h)~~~
hold\}.
\end{eqnarray*}

\par{A \textit{t-norm} $T$ is a two-place function $T:I\times I\rightarrow I$ which is associative,
 commutative, nondecreasing in each place and such that $ T(a,1)=a$, for all $a \in$ $[0, 1]$.
 A {\it triangle function} is a binary operation on $\Delta^+$, namely a function
  $\tau:\Delta^+\times\Delta^+\rightarrow\Delta^+$
that is associative, commutative, nondecreasing and which has
$\mathcal H_0$ as unit. That is, for all $F,G,H\in\Delta^{+}$, we
have}
\begin{eqnarray*} \tau(\tau(F,G),H)&=&\tau(F,\tau(G,H)),\\
\tau(F,G)&=&\tau(G,F),\\
F\leq G &\Longrightarrow & \tau(F,H)\leq\tau(G,H),\\
\tau(F,\varepsilon_{0})&=&F.\end{eqnarray*} Continuity of a triangle
functions means continuity with respect to the topology of weak
convergence in $\Delta^{+}$. Typical continuous triangle functions
are
$$\tau_{T}(F,G)(x)=\sup_{s+t=x}T(F(s),G(t)),$$ and
$$\tau_{T^{*}}(F,G)(x)=\inf_{s+t=x}T^{*}(F(s),G(t)).$$Here $T$ is
a continuous t-norm, i.e. a continuous binary operation on $[0,1]$
that is commutative, associative, nondecreasing in each variable and
has 1 as identity, and $T^{*}$ is a continuous t-conorm, namely a
continuous binary operation on $[0,1]$ which is related to the
continuous t-norm $T$ through $$ T^{*}(x,y)=1-T(1-x,1-y).$$

\begin{definition}A \textit{probabilistic normed space} is a
quadruple $(V,\nu,\tau,\tau^{*})$, where $V$ is a real vector space,
$\tau$ and $\tau^{*}$ are continuous triangle functions, and $\nu$
is a mapping from $V$ into $\Delta^{+}$ such that, for all $p,q$ in
$V$, the following conditions hold:

\,  (PN1) $\nu_{p}=\varepsilon_{0}$ if, and only if, $p=\theta$,
where $\theta$ is the null vector in $V$;

\,  (PN2) $\nu_{-p}=\nu_{p},$ for each $p\in V$;

\,  (PN3) $\nu_{p+q}\geq \tau(\nu_{p},\nu_{q}),$ for all $p,q\in V$;

\,  (PN4) $\nu_{p}\leq \tau^{*}(\nu_{\alpha p},\nu_{(1-\alpha)p})$,
for all $\alpha$ in $[0,1]$.\end{definition} If the inequality (PN4)
is replaced by the equality $$\nu_{p}=\tau_{M}(\nu_{\alpha
p},\nu_{(1-\alpha)p}),$$ then the probabilistic normed space is
called \textit{\v{S}erstnev space }and, as a consequence, a
condition stronger than (PN2) holds, namely $$ \nu_{\lambda
p}(x)=\nu_{p}(\frac{x}{|\lambda|}),$$ for all $p\in V$, $\lambda
\neq 0$ and $x\in \mathbb{R}$.

\begin{definition}
 \cite{cho} Let  $L$ be a real linear space with dim $L\geq n$ and
 $\|.,\ldots,.\|: L^n\rightarrow\mathbb{R}$ a function. Then
 $(L,\|.,\ldots,.\|)$ is called a linear $n$-normed space if

 \, (n-N1) $\|x_1,x_2,\ldots,x_n\|= 0$ if, and only if $x_1,x_2,\ldots,x_n$
are linearly dependents;

 \, (n-N2) $\|x_1,x_2,\ldots,x_n\|=\|x_{i_1},x_{i_2},\ldots,x_{i_n}\|$ for
 every permutation \\$(i_1,\ldots,i_n)$ of $(1,\ldots,n)$;

 \, (n-N3) $\|\alpha x_1,x_2,\ldots,x_n\|=|\alpha|\|x_1,x_2,\ldots,x_n\|$;

 \, (n-N4) $\|x+y,x_2,\ldots,x_n\|\leq
 \|x,x_2,\ldots,x_n\|+\|y,x_2,\ldots,x_n\|$,\\
for all $\alpha \in \mathbb{R}$ and all $x,y,x_1,\ldots, x_n \in L$.
The function $\|.,\ldots,.\|$ is called the $n$-norm on L.
 \end{definition}

\section{Main results}

\begin{definition}
Let $L$ be a linear space with dim $L\geq n$, $\tau$ a triangle
function, and let $\mathcal F$ be a mapping from $L^n$ into
$\Delta^+$. If the following conditions are satisfied :

 \,  (Pn-N1) $F_{x_1,x_2,\ldots,x_n}=\mathcal H_0$ if $x_1,x_2,\ldots,x_n$  are linearly dependent;

 \,  (Pn-N2) $F_{x_1,x_2,\ldots,x_n}\neq \mathcal H_0$  if $x_1,x_2,\ldots,x_n$  are linearly independent;

 \,  (Pn-N3) $F_{x_1,x_2,\ldots,x_n}=F_{x_{j_1},x_{j_2},\ldots,x_{j_n}}$, for every permutation $(j_1,\ldots,j_n)$;

 \,  (Pn-N4) $F_{\alpha x_1,x_2,\ldots,x_n}(t)=F_{x_1,x_2,\ldots,x_n}(\frac{t}{|\alpha|}),$ for every $t>0, \alpha\neq0$
 and $x_1,x_2,\ldots,x_n\in L$;

 \,  (Pn-N5) $F_{x+y,x_2,\ldots,x_n}\geq \tau(F_{x,x_2,\ldots,x_n},F_{y,x_2,\ldots,x_n})$, whenever
 $x,y,$\\$x_2,\ldots,x_n\in L$,\\
then $\mathcal F$ is called a {\it  probabilistic n-norm} on $L$ and
the triple ($L, \mathcal F, \tau$)~is called a {\it  probabilistic
$n$-normed space (Briefly $P$-$nN$space)}.
\end{definition}
Suppose that $L$ be a vector space with dimension $d$, where $2\leq
d< \infty$, unless otherwise stated. Fix $\{u_1,\ldots,u_d \}$ to be
a basis for $L$. Then we have the following:
\begin{lemma}\label{lem1}
Let $(L, \mathcal F, \tau)$ with continuous $\tau$ be a
$P$-$nN$space. A sequence $\{x_m\}$ in $L$ is convergent to $x$ in
$L$ if, and only if
\begin{eqnarray*}
\lim_{m\rightarrow \infty} F_{x_m-x,y_2,\ldots,y_{n-1},u_i}=\mathcal
H_0,~~~~~ for~~~ every~~~~y_2,\ldots,y_{n-1}\in L,i=1,\ldots,d.
\end{eqnarray*}
\end{lemma}
Following lemma \ref{lem1}, we have:
\begin{lemma}\label{lem2}
Let $(L, \mathcal F, \tau)$ with continuous $\tau$ be a $P$-$nN$
space. A sequence $\{x_m\}$ in $L$ is convergent to $x$ in $L$ if,
and only if
\begin{eqnarray*}
\lim_{m\rightarrow \infty}
max\{F_{x_m-x,y_2,\ldots,y_{n-1},u_i}:y_2,\ldots,y_{n-1}\in
L,i=1,\ldots,d\}=\mathcal H_0.
\end{eqnarray*}
\end{lemma}
 Now with respect to the base is $\{u_1,\ldots,u_d\}$ we can define a
norm on $L$ which we shall denote it by $\mathcal F^\infty _x$ as
follows:
\begin{eqnarray*}
\mathcal F^\infty _x:=
max\{F_{x,y_2,\ldots,y_{n-1},u_i}:y_2,\ldots,y_{n-1}\in
L,i=1,\ldots,d\}.
\end{eqnarray*}
In fact it's not difficult to see:\\
$1)$ $F^\infty _x=\mathcal H_0$ if, and only if $x=0$,\\
$2)$ $F^\infty _{\alpha x}(t)=F^\infty _x(\frac{t}{|\alpha|})$,\\
$3)$ $F^\infty _{x+y} \geq \tau (F^\infty _x,F^\infty _y)$, for $x,y
\in L$, $\alpha \neq0$, $t>0$.

Note that if we choose another basis for $L$, say
$\{v_1,\ldots,v_d\}$, and define the norm $\mathcal F^\infty$ with
respect to it, then the resulting norm will be equivlent to the one
defined with respect to $\{u_1,\ldots,u_d\}$.\\
Using the $\mathcal F^\infty$ in lemma $\ref{lem2}$ we can write it
as follows:
\begin{lemma}\label{lem3}
A sequence $\{x_m\}$ in $L$ is convergent to $x$ in $L$ if, and only
if $~~\lim_{m\rightarrow \infty} F^\infty_{x_m-x}=\mathcal H_0$.
\end{lemma}
Also associated to the derived norm $\mathcal F^\infty$, we can
define the open balls $B_{\{u_1,\ldots,u_d\}}(x,t)$ centered at $x$
having radius $t$ by
\begin{equation*}
B_{\{u_1,\ldots,u_d\}}(x,t):=\{y \in L:F^\infty_{x-y}(t)>1-t\}.
\end{equation*}
With respect to above balls, lemma $\ref{lem3}$ be comes:
\begin{lemma}\label{lem4}
A sequence $\{x_m\}$ in $L$ is convergent to $x$ in $L$ if, and only
if
\begin{eqnarray*}
\forall~~~\varepsilon>0~~~\exists~~~ N\in
\mathbb{N}~~~such~~~that~~~ m\geq N ~~~\Rightarrow x_m \in
B_{\{u_1,\ldots,u_d\}}(x,\varepsilon).
\end{eqnarray*}
\end{lemma}
Summarizing all these results, we have:
\begin{theorem}
any finite dimensional $P$-$nN$ space is a $P$-$N$ space and it's
topology agrees with that generated by the derived norm $\mathcal
F^\infty$.
\end{theorem}
\begin{definition}
A sequence $\{x_m \}$ in $P$-$nN$ space ($L,\mathcal F,\tau$) is
cauchy sequence if $\lim_{m,r\rightarrow
\infty}F_{x_m-x_r,y_2,\ldots,y_{n-1},y}=\mathcal H_0$, for every
$y\in L$.
\end{definition}
\begin{definition}
The probabilistic $n$-normed space ($L,\mathcal F,\tau$)is a
probabilistic $n$-Banach space if every cauchy sequence in $L$ is
convergent to some $x$ in $L$.
\end{definition}
\begin{theorem}
A probabilistic $n$-normed space ($L,\mathcal F,\tau$) is a
probabilistic $n$-Banach space if, and only if ($L,\mathcal F,\tau$)
is a probabilistic Banach space.
\begin{proof}
By lemma \ref{lem2}, the convergence in the probabilistic $n$-norm
is equivalent to that in the derived norm, it suffices to show
that $\{x_m\}$ is Cauchy with respect to the probabilistic
$n$-norm if, and only if it is Cauchy with respect to the derived
norm. But it is clear since $\{x_m\}$ is Cauchy with respect to
the probabilistic $n$-norm if, and only if\\ $\lim_{m,r\rightarrow
\infty}F_{x_m-x_r,y_2,\ldots,y_{n-1},y}=\mathcal H_0$ for every
$y_2,\ldots,y_{n-1},y\in L$, if and only if $\lim_{m,r\rightarrow
\infty}F_{x_m-x_r,y_2,\ldots,y_{n-1},u_i}=\mathcal H_0$ for every
$i=1,\ldots,d$, if, and only if
\begin{equation*}
\lim_{m,r\to \infty}F^{\infty}_{x_m-x_r}=\mathcal H_0
\end{equation*}
if, and only if $\{x_m\}$ is Cauchy with respect to derived norm.
\end{proof}
\end{theorem}


\end{document}